\documentclass[11pt]{article}
\usepackage[dvips]{graphicx}
\usepackage{amssymb,epsf,amsfonts,amsmath,amsxtra}

\newtheorem{thm}{Theorem}
\newtheorem{lem}[thm]{Lemma}

\newtheorem{prop}[thm]{Proposition}

\newcommand{\proof}{\noindent{\bf Proof.\ }}

\newcommand{\qed}{\begin{flushright} $\Box$ \end{flushright}}

\newcommand{\cM}{{\cal M}}
\newcommand{\cMtwo}{{\cM}^2}

\newcommand{\Spec}{{\cal S}}

\newcommand{\diagramSized}[2]{\begin{center} \includegraphics[width=#2]{#1} \end{center}}

\begin{document}

\title{On graphs having maximal independent sets\\ of exactly $t$ distinct cardinalities}

\author{$^1$Bert L. Hartnell \, and $^2$Douglas F. Rall\thanks{Corresponding author: e-mail:
        doug.rall@furman.edu }
\thanks{Research supported in part by the
Wylie Enrichment Fund of Furman University.} \\
    \\
    $^1$Department of Mathematics\\ \& Computing Science\\
        Saint Mary's University\\
        Halifax, Nova Scotia, Canada\\
    \\
    $^2$Department of Mathematics \\
    Furman University\\
    Greenville, SC 29613 USA\\
    }
    \date{}
    \maketitle

\begin{abstract}
For a given positive integer $t$ we consider graphs having maximal independent sets of precisely $t$
distinct cardinalities and restrict our attention to those that have
no vertices of degree one.  In the situation when $t$ is four or larger
and the length of the shortest cycle is at least $6t-6$, we completely
characterize such graphs.
\end{abstract}

{\small {\bf Keywords:} maximal independent set, girth, cycle} \\
\indent {\small {\bf AMS subject classification: 05C69, 05C38}}

\section{Introduction}
A well-covered graph (Plummer~\cite{mdp1970}) is one in which every maximal independent
set of vertices is of one cardinality and is hence a maximum independent set.  Finbow,
Hartnell and Whitehead~\cite{fhw1994} defined the class $\cM_t$ to consist
of those graphs which have exactly $t$ different sizes of maximal independent sets.
Finbow, Hartnell and Nowakowski~\cite{fhn1993} proved that the well-covered graphs
(the $\cM_1$ collection) of girth (the length of a shortest cycle) 6 or more,
with the exceptions of $K_1$ and $C_7$,
have the property that every vertex has degree one or has exactly one vertex of
degree one in its neighborhood.  Thus, $C_7$ is the unique graph in $\cM_1$ with
girth at least 6 that has minimum degree at least two.  The graphs in $\cM_2$
of girth 8 or more have also been characterized (\cite{fhw1994}).  There are
precisely five graphs in $\cM_2$ of girth at least 8 and minimum degree 2 or more,
namely the cycles $C_8, C_9, C_{10}, C_{11}$ and $C_{13}$.  This implies there
are no $\cM_1$ graphs of girth at least 8 with minimum degree 2 or more and no
$\cM_2$ graphs of girth 14 or more and having minimum degree at least 2.  For
related work on the class $\cM_t$ see \cite{bh1998} and \cite{bh2009}.

In this paper we investigate the graphs in $\cM_t$ that have minimum degree
at least 2 and higher girth and establish that the characterization of these
in $\cM_1$ and $\cM_2$ is part of a general pattern.  In particular, for $t \ge 3$ we show
that among graphs with minimum degree at least 2, $\cM_t$ does not contain
a graph of girth at least $6t+2$ and that $C_{6t-4}, C_{6t-3}, C_{6t-2},
C_{6t-1}$ and $C_{6t+1}$ are the only exceptions for girth at least $6t-4$.
Furthermore, if $t \ge 4$, then these cycles along with $C_{6t-6}$ are the
only graphs in $\cM_t$ that have minimum degree at least 2 and girth at least
$6t-6$.

Let $G$ be a finite simple graph.  A vertex of degree 1 is called
a {\it leaf} and any vertex that is adjacent to a leaf is called a
{\it support vertex}.  If $C$ is a cycle in a graph $G$ and $u$ and $v$
belong to $C$, we let $uCv$ denote the shorter of the two $u,v$-paths that
are part of $C$.   For $A \subseteq V(G)$ and $u$ a vertex in $G$, $d(u,A)$ will
denote the length of a shortest path in $G$ from $u$ to a vertex of $A$.
We will use $\cM(G)$ to denote
the collection of all maximal independent sets of $G$ and we define the {\it independence
spectrum} ({\it spectrum} for short) of $G$ to be the set
$\Spec(G)=\{\,|I|\,:\,I\in \cM(G)\}$.  The class $\cM_t$ consists of those
graphs $G$ for which $|\Spec(G)|=t$.  The spectrum is not necessarily a set
of consecutive positive integers (e.g., $\Spec(K_{2,4,5})=\{2,4,5\}$), but for
paths and cycles it is.  We denote the set of positive integers between $p$ and
$q$ inclusive by $[p,q]$.  The following proposition is easy to establish.

\begin{prop}\label{pathcycle}
For each positive integer $n$ at least 3,
\[\Spec(C_n)=[\lceil{n/3}\rceil,\lfloor{n/2}\rfloor]\quad \mbox{and} \quad \Spec(P_n)
=[\lceil{n/3}\rceil,\lceil{n/2}\rceil]\,.\]  Hence, $C_n \in \cM_t$ and $P_n \in \cM_s$
where $t=\lfloor{n/2}\rfloor - \lceil{n/3}\rceil +1$ and $s= \lceil{n/2}\rceil -
\lceil{n/3}\rceil +1$.
\end{prop}

The following lemma from \cite{fhw1994} will be used throughout---often without mention.

\begin{lem} \label{leftover} {\rm \cite{fhw1994}}
If the graph $G$ belongs to $\cM_t$ and $I$ is an independent set of $G$,
then for every component $C$ of $G-N[I]$ there exists $k\le t$ such that
$C \in \cM_k$.  In addition, $G-N[I] \in \cM_r$ for some $r \le t$.
\end{lem}

Lemma~\ref{leftover} will most often be used in the following way.  We will
find an independent set $I$ in a graph $G$ and demonstrate that $G-N[I]$
has a component that is in the class $\cM_s$ for some $s>t$ and conclude
that $G \not\in\cM_t$.  The following lemma will be used in that
context with Lemma~\ref{leftover}.

\begin{lem} \label{cycleplus1}
If a cycle $C$ is in $\cM_t$ and a new vertex is added as a leaf adjacent
to a single vertex of $C$, then the resulting graph belongs to $\cM_{t+1}$.
\end{lem}
\proof  Assume $\Spec(C)=[k,k+t-1]$.  Let $H$ be the graph formed by
adding a leaf $x$ adjacent to $y$.  Let $u$ and $v$ be the neighbors of $y$
on $C$.  Note that $\{ I \in \cM(H)\,:\, y \in I\}=\{ J \in \cM(C)\,:\, y \in J\}$,
and because of the symmetry of the cycle, $\Spec(C)=\{|J|\,:\, J \in \cM(C), y \in J\}$.
Also, $\{ I \in \cM(H)\,:\, u \in I\}=\{ J \cup\{x\}\,:\,  J \in \cM(C), u \in J\}$.
This shows that $[k,k+t] \subseteq \Spec(H)$.  If $H$ has a maximal independent set $A$
of size less than $k$, then $x\in A$ and neither $u$ nor $v$ is in $A$, for otherwise
$A \cap C$ is a maximal independent set in $C$ of cardinality less than $k$.  But now
$A'=(A-\{x\})\cup \{y\} \in \cM(C)$ and $|A'|<k$, a contradiction.  Therefore,
$\Spec(H)=[k,k+t]$. We conclude that $H \in \cM_{t+1}$. \qed

In the class of graphs with leaves there is no connection between girth and
the size of the spectrum.  This can be seen by the following general construction.
Let $t\ge 2$ and $g\ge 3$ be integers.  Let $H$ be the graph formed by adding a single
leaf adjacent to each vertex of a cycle of order $g$.  For a single vertex $x$ on the
cycle attach a path $v_1, v_2, \ldots, v_{2t-3}$ to $H$ by making $x$ and $v_1$ adjacent.
Then add two leaves adjacent to $v_i$ if $i$ is odd, and add one leaf adjacent to
$v_j$ if $j$ is even. The resulting graph of order $2g+5t-7$ has girth $g$ and belongs
to the class $\cM_t$. (The spectrum of this graph is $[g+2t-3, g+3t-4]$.)  For this
reason we will henceforth consider only graphs having minimum degree at least 2.  For
ease of reference we denote the class of graphs that are in $\cM_t$ and have
no leaves (i.e., minimum degree at least 2) by $\cMtwo_t$.  Note that $\cMtwo_t \subseteq \cM_t$.
In the course of several of our proofs we will show that some given graph is
not in $\cMtwo_t$ by demonstrating it does not belong to $\cM_t$.

The remainder of this paper is devoted to verifying the entries in the following
table.
\noindent
\begin{table}[h]
\begin{center}
\begin{tabular}{ | c|c|c|c|c|c|c|c|c|c|}
\hline
&\multicolumn{9}{c |} {\emph{girth}} \\  \hline
     & $6t-6$ & $6t-5$ & $6t-4$ & $6t-3$ & $6t-2$ & $6t-1$ & $6t$ & $6t+1$ & $\ge 6t+2$  \\ \hline
\hline
$t=1$ &   &    &   &  $\Delta$ & $\Delta$ & $\Delta$ & $\emptyset$ & $C_7$  & $\emptyset$ \\ \hline
$t=2$ & $\Delta$  & $\Delta$   & $C_8$  &  $C_9$ & $C_{10}$ & $C_{11}$ & $\emptyset$ & $C_{13}$  & $\emptyset$ \\ \hline

$t=3$ & $C_{12}$  & $\Delta$   & $C_{14}$  &  $C_{15}$ & $C_{16}$ & $C_{17}$ & $\emptyset$ & $C_{19}$  & $\emptyset$ \\ \hline

$t=4$ & $C_{18}$  & $\emptyset$   & $C_{20}$  &  $C_{21}$ & $C_{22}$ & $C_{23}$ & $\emptyset$ & $C_{25}$  & $\emptyset$ \\ \hline
$t\ge 5$ & $C_{6t-6}$  & $\emptyset$   & $C_{6t-4}$  &  $C_{6t-3}$ & $C_{6t-2}$ & $C_{6t-1}$ & $\emptyset$ & $C_{6t+1}$  & $\emptyset$ \\ \hline
\end{tabular}
\caption{Graphs of given girth in $\cMtwo_t$} \label{summarytable}
\end{center}
\end{table}

The entry for a given girth (written as a function of $t$) and a given
value of $t$ should be interpreted as follows.  If a specific graph is given, then
this is the unique graph of that girth that belongs to $\cMtwo_t$. For example,
$C_{15}$ is the only graph of girth 15 in $\cMtwo_3$. If $\emptyset$
appears, then there are no graphs of that girth in $\cMtwo_t$.  When the entry
is $\Delta$, then it is known that $\cMtwo_t$ contains at least one graph of that
girth (and it is not just a cycle).  Some of these type of entries have been verified in previous papers.  For
example, see \cite{fhn1993} and \cite{fhw1994} for $\cMtwo_1$ and $\cMtwo_2$, respectively.

\section{Establishing Table Entries}

We begin by showing that for a given positive integer $t$ the only
graphs in $\cM_t$ with large enough girth must have leaves.  The next result was
proved for well-covered graphs ($t=1$) in \cite{fh1983}.  Proposition \ref{pathcycle}
shows it is sharp in terms of girth.

\begin{thm} \label{leaves}
Let $t$ be a positive integer.  If $g(G)\ge 6t+2$ and $\delta(G) \ge 2$, then $G \in \cM_r(G)$ for some $r>t$.
\end{thm}
\proof  Assume $t\ge 2$. Let $G$ have girth at least $6t+2$ and minimum degree at least two.  We
will show that $G$ has maximal independent sets of at least $t+1$ different sizes.  Choose a cycle
$C=v_1,v_2,\ldots,v_s$ of minimum length in $G$.

Assume first that $s \ge 6t+4$ and let $P$ denote the path $v_3,v_4,\ldots,
v_{6t+1}$.  Since $\delta(G) \ge 2$ and $g(G)=s$, each vertex $u\not\in C$
that is adjacent to a vertex of $P$ has another neighbor $u'$ that does not
belong to $P$ and is not adjacent to any vertex of $P$.  Choose one such neighbor $u'$
for each $u$ and let $J$ denote the set of these neighbors.  By the girth restriction it follows that the
set $I=J\cup\{v_1,v_{6t+3}\}$ is independent. (If $s=6t+2$, then proceed as above except let
$I=J\cup \{v_1\}$.)  However, $P$ is a component of
$G-N[I]$ and by Proposition \ref{pathcycle}, $P \in \cM_{t+1}$.  Similar to the
proof of Lemma~\ref{leftover} this
implies that $G$ has maximal independent sets of at
least $t+1$ different sizes.

If $s=6t+3$, let $P$ be the path $v_3,v_4,\ldots, v_{6t+2}$.  The set $J$ is
chosen as before, and now $G-N[J \cup \{v_1\}]$ has the path $P$ of order $6t$
as a component.  By Proposition~\ref{pathcycle} it once again follows that
$G$ has at least $t+1$ distinct sizes of maximal independent sets.  \qed
\medskip
For any positive integer $t$ it follows from Proposition~\ref{pathcycle}
that $C_{6t+1} \in \cM_t$.   In \cite{fhn1993} it was shown that $C_7$ is the
only well-covered graph of girth 7 and minimum degree 2 or more.  The following
theorem shows the similar result is true for larger values of $t$.

\begin{thm}\label{girth6tand6t+1}
Let $t\ge 2$ be an integer.  The cycle $C_{6t+1}$ is the only graph of girth
$6t+1$ in $\cMtwo_t$, and $\cMtwo_t$ contains no graphs of girth $6t$.
\end{thm}
\proof By Proposition~\ref{pathcycle} the cycle of order $6t+1$ belongs to
$\cMtwo_t$.  Suppose $G$ is a graph not isomorphic to $C_{6t+1}$ such that
$g(G)=6t+1$ and $\delta(G) \ge 2$.  Then $G$ has an induced cycle $C$ of order $6t+1$,
and $C$ has a vertex $w$ of degree at least 3.
Since $g(G)=6t+1$ and $\delta(G)\ge 2$ we can find an induced path $w,a,b,c$,
such that none of $a, b$ or $c$ belongs to $C$.  Let $X =\{ u \in V(G)\, :\, d(u,C)=2 \}-N(a)$
and let $Y=\{ u \in V(G)\,:\, d(u,a)=2, d(u,w)=3\}$. For any two vertices on $C$
there is a path using part of $C$ of length at most $3t$ joining them.
Since $g(G) \ge 13$ it follows that $Y$ is independent. Suppose
two vertices $x_1,x_2 \in X$ are adjacent.  Let $x_1,v_1,w_1$ and $x_2,v_2,w_2$
be paths in $G$ with $w_1$ and $w_2$ on the cycle $C$.  Then the cycle
$x_1,v_1,w_1Cw_2,v_2,x_2,x_1$ has length at most $3t+5$.  But then $3t+5\ge 6t+1$,
which implies that $t=1$, a contradiction.    Finally, if a vertex in $X$ is
adjacent to a vertex in $Y$, then a similar argument shows that $G$ has a
cycle of length at most $3t+6$ which also leads to a contradiction.

Therefore, $X \cup Y$ is an independent set.  One of the components of the
graph $G-N[X \cup Y]$ is the cycle $C$ with a single leaf $a$ attached at the
support vertex $w$.  By Lemma~\ref{cycleplus1} this component is in $\cM_{t+1}$.
An application of Lemma~\ref{leftover} then shows that $G \not\in\cMtwo_t$.

Now let $G$ be a graph of girth $6t$, and as above find an induced cycle $C$
of length $6t$.  This time let  $X =\{ u \in V(G)\, :\, d(u,C)=2 \}$.  This set
is independent unless there is a cycle of the form $x_1,v_1,w_1Cw_2,v_2,x_2,x_1$
that has length at most $3t+5$.  But this means $3t+5 \ge 6t$ contradicting our
assumption that $t\ge 2$.  Hence $X$ is independent.  The cycle $C$ is one of the
components of $G-N[X]$.  Since $C_{6t} \in \cM_{t+1}$, Lemma~\ref{leftover} implies
that $G \not\in \cMtwo_t$.    \qed

By following a line of reasoning similar to the first part of the proof of
Theorem~\ref{girth6tand6t+1} one can prove the following result.  The proof
is omitted.  As noted earlier, Theorem~\ref{fourgirths} also holds for $t=2$.  See \cite{fhw1994}.

\begin{thm}\label{fourgirths}
Let $t \ge 3$ be a positive integer.  For each integer $n$ such that
$6t-4 \le n \le 6t-1$, the cycle $C_n$ is the unique graph
of girth $n$ that belongs to $\cMtwo_t$.
\end{thm}

We now establish the uniqueness (for $t \ge 3$) of the table entry corresponding to those
graphs with no leaves whose shortest cycle has length $6t-6$ and which have maximal
independent sets of exactly $t$ distinct cardinalities.

\begin{thm}\label{girth6t+1}
For each integer $t \ge 3$, the cycle $C_{6t-6}$ is the only graph of girth
$6t-6$ that belongs to $\cMtwo_t$.
\end{thm}
\proof  The cycle of order $6t-6$ is in $\cMtwo_t$ by Proposition~\ref{pathcycle}.
Suppose that $G$ is a graph of girth $6t-6$  with no leaves.  If $G$ is not $C_{6t-6}$,
then  we can find an induced cycle $C$ of length $6t-6$ in $G$ with $w,a,b,c$, $X$ and $Y$
defined as in the proof of Theorem~\ref{girth6tand6t+1}.  The set $Y$ is independent because
$g(G)\ge 12$, and $X$ is independent since $t \ge 3$.  If some vertex of $X$ is adjacent to a
vertex of $Y$, then $G$ contains a cycle of length at most $3t-3+6$.  It follows that
$3t+3 \ge g(G)= 6t-6$, or equivalently $t \le 3$.

If the set $X \cup Y$ is independent, then
$G-N[X \cup Y]$ has a component isomorphic to a cycle of length $6t-6$ with a single leaf
attached at $w$.  By Lemma~\ref{cycleplus1} this component is in $\cM_{t+1}$
and so it follows from Lemma~\ref{leftover} that $G \not\in\cM_t$.

\begin{figure}[h]
\diagramSized{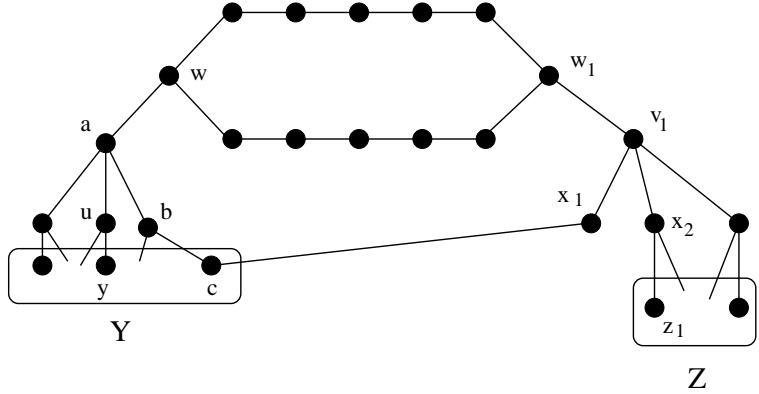}{10 truecm}
\caption{Part of $G$}
\label{Fig1}
\end{figure}

Thus we may assume that $t=3$ and that $X \cup Y$ is not independent.  Without loss of
generality we may assume that $c$ from $Y$ is adjacent to $x_1$ such that $x_1\in X$
and $x_1,v_1,w_1$ is a path where $w_1$ is on the cycle $C$.  See Figure~\ref{Fig1}.
By using the fact that $C$
has length 12 and $g(G)=12$ we infer that the length of $wCw_1$ is 6.  Let $X'=X-N(v_1)$ and
let $Z=\{u\,:\,d(u,v_1)=2, d(u,w_1)=3, ux_1 \not\in E(G)\}$.  It is clear that $Z$ is independent.

As above, if a vertex of $Z$ is adjacent to a vertex $h$ of $X'$, then if $d(h,w)>2$ a cycle
of length at most 11 is present and if $d(h,w)=2$ then $G$ contains a cycle of length 10,
contradicting $g(G)=12$.  Suppose $z_1 \in Y\cap Z$, say $z_1=y$ as in Figure~\ref{Fig1}.
Then $z_1 \not =c$, and $a,b,c,x_1,v_1,x_2,z_1,u,a$ is a cycle, contradicting the girth
assumption.  Similarly, since $G$ has no cycles of length 9, it follows that $Z \cup Y$ is
independent.

The set $X' \cup Y \cup Z$ is independent, and one of the components of the graph
$G-N[X' \cup Y \cup Z]$ is the cycle $C$ with a single leaf attached at vertices
$w$ and $w_1$.  But this component has spectrum $\{4,5,6,7,8\}$ from which it follows
that $G\not\in\cM_3$. \qed

We now show that when $t \ge 4$ there is a ``gap'' at girth $6t-5$ among the leafless graphs.  That is,
if $G$ has minimum degree at least 2 and the shortest cycle of $G$ has order $6t-5$,
then $G$ does not belong to $\cM_t$.

\begin{figure}[htb]
\diagramSized{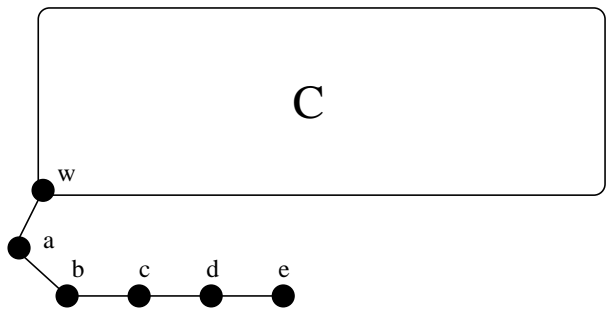}{8 truecm}
\caption{The cycle $C$ with attachments}
\label{Fig2}
\end{figure}

\begin{thm} \label{lowergap}
For each integer $t$ at least 4, the class $\cMtwo_t$ contains no graphs of girth $6t-5$.
\end{thm}
\proof  First observe that $C_{6t-5}\in \cM_{t-1}$.  Our approach will be similar as
that pursued in earlier proofs, except that we will be attempting to isolate a
cycle of length $6t-5$ with a path of order 5 attached as in Figure~\ref{Fig2}.  It
is easy to check, using either $\{a,c,e\}$ or $\{a,d\}$ together with all
possible maximal independent sets of a path of order $6t-6$, that this component has
spectrum $[2t,3t]$ and hence belongs to $\cM_{t+1}$.  This in turn implies via Lemma~\ref{leftover}
that $G \not\in\cMtwo_t$.

Suppose that $G$ has girth $6t-5$ and has minimum degree at least 2.  Let $C$ be an induced cycle
of length $6t-5$ in $G$.  There must exist a vertex $w$ on $C$ having degree at least 3.  For any
two vertices on $C$ there is a path on $C$ joining them whose length is at most $3t-3$.
Because of the girth and minimum degree assumptions on $G$ we can find a path
$w,a,b,c,d,e$ as in Figure~\ref{Fig2}.  Let $A=\{a,b,c,d,e\}$.  Let $X=\{u\,:\,d(u,C)=2\}-N(a)$
and let $Y=\{u\,:\,u \not\in C, d(u,A)=2, d(u,w)\ge 2\}$.

As in previous proofs it is straightforward to show that $X$ is independent.  Since
$g(G)=6t-5\ge 19$ no pair of vertices in $Y$ can be adjacent.  Suppose first that $X\cup Y$
is independent.  The graph in Figure~\ref{Fig2} is a component of $G-N[X\cup Y]$.  As
remarked at the outset, this shows that $G \not\in\cMtwo_t$.  We note that for $t\ge 5$,
the girth restriction ensures that $X\cup Y$ is independent.

Now consider $t=4$.  Thus $C$ is of length 19.  Let $s_1$ and $s_2$ be the adjacent vertices
on $C$ that are at distance 9 from $w$.  If both $s_1$ and $s_2$ are of degree two, then
$X \cup Y$ is independent or else a cycle of length 18 would exist in $G$. Assume then without
loss of generality that $s_1$ has a neighbor $r$ that is not on $C$. Let $U=N(r)-\{s_1\}$.
For each $u_i \in U$ choose a vertex $v_i \in N(u_i)-\{r\}$, and set
$V=\{v_i\,:\, u_i\in U\}$.  Similarly, let $B=N(a)-\{w\}$.  For each $b_i \in B$ choose a
vertex $c_i \in N(b_i)-\{a\}$, and set $D=\{c_i\,:\, b_i\in B\}$.  Since
$g(G)=19$ the set $V \cup D \cup (X-U)$ is independent, and one of the components
of $G-N[V \cup D \cup (X-U)]$ is a cycle of order 19 with a single
leaf $a$ adjacent to $w$ and a single leaf $r$ adjacent to $s_1$.  This component
belongs to $\cM_5$ which proves that $G \not\in\cMtwo_4$ and establishes the theorem. \qed

\section{Concluding Remarks}

We have shown that for a positive integer $t\ge 4$ and for each possible value of girth
at least $6t-6$, the class $\cMtwo_t$ either contains exactly one graph of that girth (the cycle)
 or contains no graphs of that girth.  It is interesting to note that as $t$ grows there
 is an ever increasing gap---in terms of girth---between the unique graph of girth $6t-6$ in
 $\cMtwo_t$ and ones of smaller girth.  For instance, we can show that $\cMtwo_{31}$ contains
 no graphs of girth $r$ for $131 \le r \le 179$.  Hence the cycles $C_{180}, C_{182}, C_{183},
 C_{184}, C_{185}$ and $C_{187}$ are the only leafless members of $\cM_{31}$ that have girth
 at least 131.  Thus the six cycles are quite special in $\cMtwo_t$.

\end{document}